\documentclass[dvips,color,xnumber,nohams]{amsart}
\usepackage{cras,math,curves,eepic,epsfig,t1enc}
\newcommand\MAT[4]{\begin{bmatrix}#1 & #2\\#3 & #4\end{bmatrix}}
\renewcommand\gf{{\mathcal G}}
\newtheorem*{defnn}{Definition}

\begin{document}
\issueinfo{33?}{?}{?}{2000}
\title[Spectra of graphs related to fractal groups]{Spectra of
  non-commutative dynamical systems and\\
  graphs related to fractal groups}
\author{Laurent BARTHOLDI}
\address{Section de Math\'ematiques, Universit\'e de Gen\`eve, CP 240,
  1211 Gen\`eve 24, Suisse}
\email{Laurent.Bartholdi@math.unige.ch}
\author{Rostislav I. GRIGORCHUK}
\address{Steklov Institute of Mathematics, Gubkina 8, Moscow 117966, Russia}
\email{grigorch@mi.ras.ru}
\date{2 juin 2000}
\commby{Jean-Pierre \textsc{Serre}}
\seriesinfo{Théorie des groupes}{Group Theory}
\begin{abstract}
  We study spectra of noncommutative dynamical systems,
  representations of fractal groups, and regular graphs. We explicitly
  compute these spectra for five examples of groups acting on rooted
  trees, and in three cases obtain totally disconnected sets.
\end{abstract}
\ftitle{Spectres de systèmes dynamiques non-commutatifs et
  de graphes associés à des groupes fractals}
\begin{fabstract}
  Nous étudions les spectres de systèmes dynamiques non commutatifs,
  de représentations de groupes fractals, et de graphes réguliers.
  Nous calculons ces spectres pour cinq exemples de groupes agissant
  sur des arbres enracinés, et dans trois cas obtenons des ensembles
  complètement déconnectés.
\end{fabstract}
\maketitle

%%%%%%%%%%%%%%%%%%%%%%%%%%%%%%%%%%%%%%%%%%%%%%%%%%%%%%%%%%%%%%%%
\begin{frenchversion}
  Le but de cette note est d'étudier le spectre de systèmes dynamiques
  non commutatifs, i.e.\ de systèmes engendrés par plusieurs
  transformations qui ne commutent pas nécessairement. Elle résume les
  résultats de~\cite{bartholdi-g:spectrum}.
  
  Soit $G$ un groupe engendré par un ensemble symétrique fini $S$.  Un
  \emph{système dynamique}, noté $(S,X,\mu)$, est une action de $G$
  sur un espace mesuré $X$ préservant la classe $[\mu]$ de la mesure
  $\mu$ de $X$. Il induit naturellement une représentation unitaire
  $\pi$ de $G$ dans $L^2(X,\mu)$ donnée par $(\pi(g)f)(x) =
  \sqrt{\mathfrak g(x)}f(g^{-1}x)$, où $\mathfrak g(x)$ est la dérivée
  de Radon-Nikod\'ym $dg\mu(x)/d\mu(x)$. Le \emdef{spectre} de
  $(S,X,\mu)$ est le spectre de l'opérateur de type Hecke
  \[H_\pi = \frac1{|S|}\sum_{s\in S}\pi(s)\in\mathcal B(L^2(X,\mu)).\]
  Plus généralement, le spectre d'une représentation unitaire
  $\pi:G\to\uhilb$ d'un groupe avec système générateur fixé est le
  spectre de l'opérateur $H_\pi$ comme ci-dessus.
  
  On considère, pour cinq exemples de groupes donnés par leur action
  sur un arbre régulier enraciné $\tree$, leur représentation $\pi$
  dans $L^2(\partial\tree,\mu)$ où $\mu$ est la mesure de Bernoulli
  uniforme. Cette représentation s'approxime par les représentations
  $\pi_n$ sur les sommets à distance $n$ de la racine de $\tree$.
  
  Nous montrons que $\pi$ se décompose en une somme de représentations
  $\pi_n\ominus\pi_{n-1}$ de dimensions finies, et
  $\spec(\pi)=\overline{\bigcup_{n\ge0}\spec(\pi_n)}$.
  
  $H_\pi$ a un spectre purement ponctuel et son rayon spectral est une
  valeur propre. On peut décrire explicitement ces spectres comme
  suit. Pour $\lambda\in\R$, soit $J(\lambda)$ l'ensemble de Julia du
  polynôme quadratique $z\mapsto z^2-\lambda$ :
  $J(\lambda)=\overline{\left\{\sqrt{\lambda\pm\sqrt{\lambda\pm\sqrt{\lambda\pm\sqrt{\dots}}}}\right\}}$.

  \begin{center}
    \begin{tabular}{r|c|p{6cm}}
      Groupe & Spectre de $H_\pi$ & Description\\[1pt] \hline
      $\Gg$ & $[-\frac12,0]\cup[\frac12,1]$ &deux intervalles\rule{0pt}{12pt}\\
      $\tilde\Gg$ & $[0,1]$ &intervalle positif\\
      $\Gamma$ & $\{1,\frac14\}\cup \frac14(1\pm J(6))\}$ &union d'un
      ensemble de Cantor de mesure nulle et de points isolés\\
      $\overline\Gamma$ &
      $\{1,-\frac12,\frac14\}\cup\frac14\Big(1\pm\sqrt{\frac92\pm2J(\frac{45}{16})}\Big)$
      &ensemble de Cantor de mesure nulle\\
      $\doverline\Gamma$ & & identique à $\overline\Gamma$\\
    \end{tabular}
  \end{center}
  
  Ces calculs impliquent l'existence de graphes à croissance
  polynomiale, qui sont les graphes de Schreier de groupes de
  croissance intermédiaire, et dont le spectre de l'opérateur de
  Markov est un quelconque des ensembles ci-dessus. Un résultat
  analogue est vrai pour les systèmes dynamiques.
\end{frenchversion}
  
%%%%%%%%%%%%%%%%%%%%%%%%%%%%%%%%%%%%%%%%%%%%%%%%%%%%%%%%%%%%%%%%
\section{Introduction}
The purpose of this note is the study of spectra of noncommutative
dynamical systems (that is, systems generated by several
transformations that do not necessarily commute). We produce several
examples of computations of such spectra with an interesting
topological structure. More details appear
in~\cite{bartholdi-g:spectrum}.

In the classical case defined by a single aperiodic measure-preserving
transformation $T:X\to X$, the spectrum of the corresponding operator
$\frac12(U+U^{-1})$, with $U\in\mathcal U(L^2(X))$, is $[-1,1]$ by
Rohlin's Lemma, but in the noncommutative case the spectrum may have
gaps, and even be a Cantor set.

We also study spectra of infinite regular graphs and produce the first
example of a regular graph whose spectrum is a Cantor set.

Our dynamical systems (with associated Hecke operator) arise from
actions of fractal groups on the boundary of the regular rooted tree
on which they act. The graphs (with associated Markov operator) whose
spectra we consider are the Schreier graphs of these groups over
``parabolic subgroups''. In special cases, these graphs are
``substitutional graphs'' and have polynomial growth.

If the underlying group is amenable, then the spectrum of the the
dynamical system coincides with the spectrum of the corresponding
graph.

The computation of the above spectra is based on operator recursions
that hold for fractal groups and involves a $1$-dimensional and
$2$-dimensional classical dynamical system as an intermediate
step. This leads to the appearance of Julia sets of quadratic maps in
the description of the above spectra.

Both authors wish to thank heartily Pierre de la Harpe and Alain
Valette for their numerous comments and contributions to this note.

%%%%%%%%%%%%%%%%%%%%%%%%%%%%%%%%%%%%%%%%%%%%%%%%%%%%%%%%%%%%%%%%
\section{Spectra of Dynamical Systems and Representations}
Let $G$ be a group finitely generated by a symmetric set $S$. A
\emdef{non-commutative dynamical system}, denoted $(S,X,\mu)$, is an
action of a group $G$ (generated by $S$) on a space $X$ and preserving
the measure class $[\mu]$ of a measure $\mu$ on $X$.

Such a dynamical system gives rise to a natural unitary representation
$\pi$ of $G$ in $L^2(X,\mu)$ given by
\[(\pi(g)f)(x) = \sqrt{\mathfrak g(x)}f(g^{-1}x),\]
where $\mathfrak g(x)=dg\mu(x)/d\mu(x)$ is the Radon-Nikod\'ym
derivative.  The \emdef{spectrum} of the dynamical system $(S,X,\mu)$
is the spectrum of the Hecke type operator
\[H_\pi = \frac1{|S|}\sum_{s\in S}\pi(s)\in\mathcal B(L^2(X,\mu))\]
(this terminology comes from an analogy with Hecke operators in number
theory~\cite{serre:hecke}.)  More generally, the spectrum of a unitary
representation $\pi:G\to\uhilb$ of a group with a given finite set of
generators is the spectrum of the operator $H_\pi$ as above ---
see~\cite{harpe-r-v:sg}.

\begin{defnn}
  A \emdef{graph} is a pair $\gf=(V,E)$ of sets (the \emdef{vertices}
  and \emdef{edges}), a map $\alpha:E\to V$ (the \emdef{start} of an
  edge) and an involution $\overline\cdot:E\to E$ (the
  \emdef{inversion}). One defines then the \emdef{end} $\omega(e)$ of
  an edge by $\omega(e)=\alpha(\overline e)$.
  
  The \emdef{degree} of a vertex is $\deg(v)=|\{e\in
  E|\,\alpha(e)=v\}|$. The graph is \emdef{locally finite} if
  $\deg(v)<\infty$ for all $v\in V$, and is \emdef{regular} if
  $\deg(v)$ is constant over $V$. The graph is a \emdef{tree} if in
  every circuit $(e_1,e_2,\dots,e_n)$ of edges with
  $\omega(e_i)=\alpha(e_{i+1})$ (indices modulo $n$) there is a
  \emdef{reduction}, i.e.\ an $i$ with $\overline{e_i}=e_{i+1}$. A
  \emdef{graph morphism} is a pair of maps between the vertex and edge
  sets that commute with $\alpha$ and $\overline\cdot$.
  
  Let $\Sigma=\{1,\dots,d\}$ be a finite set of cardinality $d$. The
  \emdef{$d$-regular rooted tree} $\tree_d$ is the graph with vertex
  set $\Sigma^*$, edge set $\Sigma^*\times\Sigma\times\{\pm\}$, and
  maps $\alpha(\sigma,s,+)=\sigma$, $\alpha(\sigma,s,-)=\sigma s$ and
  $\overline{(\sigma,s,+)}=(\sigma,s,-)$.  Its \emdef{boundary} is
  $\partial\tree=\Sigma^\N$, the set of infinite sequences over
  $\Sigma$.
\end{defnn}

Suppose now that $G$ acts by automorphisms on a rooted tree $\tree$.
This action extends to a continuous action on the boundary of the
tree, which is a compact, totally disconnected space. The
\emdef{uniform measure} on $\Sigma^\N$ is the measure $\mu$ defined on
the cylinders $\sigma\Sigma^\N$ (for $\sigma\in\Sigma^*$) by
$\mu(\sigma\Sigma^*)=d^{-|\sigma|}$. It is $G$-invariant, and is the
unique invariant measure if $G$ acts transitively on each level
$\Sigma^n$ of the tree, or equivalently if $(S,\Sigma^\N,\mu)$ is
ergodic. Note that all other nondegenerate Bernoulli measures are
quasi-invariant.

For each $n\in\N$, let $\pi_n$ be the unitary representation (of
finite dimension $d^n$) of $G$ in $\ell^2(\Sigma^n)$ induced by the
action of $G$ on the $n$-th level, and let $\pi$ be the unitary
representation of $G$ in $L^2(\Sigma^\N,\mu)$.  Then $\pi_n$ is a
subrepresentation of $\pi_{n+1}$ and of $\pi$, and the spectra of
$\pi_n$ converge to that of $\pi$:
\[\spec(\pi)=\overline{\bigcup_{n\in\N}\spec(\pi_n)}.\]

%%%%%%%%%%%%%%%%%%%%%%%%%%%%%%%%%%%%%%%%%%%%%%%%%%%%%%%%%%%%%%%%
\section{Spectra of Graphs}
Let $\gf$ be a locally finite graph. The \emdef{Markov operator} of
$\gf$ is the operator $M$ on $\ell^2(V)$ given by
\[(Mf)(v)=\frac1{\deg v}\sum_{e\in E:\,\alpha(e)=v}f(\omega(e)).\]
The operator $M$ is the transition operator for the simple random walk
on $\gf$.

The spectral properties of $M$ are of great importance; for instance,
a theorem of Kesten~\cite{kesten:rwalks} (extended by Dodziuk and
others) claims that a graph $\gf$ of bounded degree is amenable if and
only if $1\in\spec(M)$.  (See more on amenability in
Section~\ref{sec:amen}.)

\begin{defnn}
  Let $G$ be a group finitely generated by a symmetric set $S$, and
  let $H$ be any subgroup. The \emdef{Schreier graph} of $G$ with
  respect to $H$ is the graph $\sch(G,H,S)$ with vertex set $G/H$ and
  edge set $S\times G/H$, and maps $\alpha(s,gH)=gH$ and
  $\overline{(s,gH)}=(s^{-1},sgH)$. It has a natural base-point $H$.

  In Subsection~\ref{subs:subs} the graphs will be labelled. This is
  simply done by assigning to each edge $(s,gH)$ the labeling $s$.
\end{defnn}
If $H=1$, we obtain the usual Cayley graph of $(G,S)$. Note that
$\sch(G,H,S)$ is an $|S|$-regular graph, but its automorphism group
does not necessarily act transitively on its vertices. Indeed,
basically any regular graph is a Schreier
graph~\cite[Theorem~5.4]{lubotzky:cayley}.

$G$ acts by left-multiplication on $G/H$, the vertex set of
$\sch(G,H,S)$. The corresponding unitary representation in
$\ell^2(G/H)$ is the \emdef{quasi-regular} representation
$\rho_{G/H}$.

Suppose now that $G$ acts on a rooted tree $\tree_d$. Fix the ray
$e=dd\dots\in\Sigma^\N$. Let $P_n$ be the stabilizer of vertex at
level $n$ of this ray, and let $P=\bigcap_{n\in\N}P_n$ be the
stabilizer of the infinite ray $e$ (it is called a \emdef{parabolic
  subgroup}). We write $M_n$ the Markov operator of $\sch(G,P_n,S)$ and
$M$ the Markov operator of $\sch(G,P,S)$; they are the Hecke
type operators associated with the action of $G$ on $\Sigma^n$ and
$\Sigma^\N$.

%%%%%%%%%%%%%%%%%%%%%%%%%%%%%%%%%%%%%%%%%%%%%%%%%%%%%%%%%%%%%%%%
\section{Fractal Groups and Substitutional Graphs}
Let $G$ be a group acting on a tree $\tree$, and let
$H=\bigcap_{s\in\Sigma}\stab_G(s)$ be the \emdef{first level
  stabilizer}. Restricting the action of $H$ to each subtree spanned
by $s\Sigma^*$ gives an embedding
\[\psi:H\to\aut(\tree)^\Sigma.\]

\begin{defnn}
  The group $G$ is \emdef{fractal} if $\psi$ is a subdirect embedding
  of $H$ in $G^\Sigma$, i.e.\ if $\psi(H)$ lies in $G^\Sigma$ and its
  projection on each factor is onto.
\end{defnn}

Let now $G$ be a group finitely generated by a symmetric set $S$, and
let $(X,x_0)$ be a pointed space on which $G$ acts. The \emdef{growth}
of $X$ is the function $\gamma:\N\to\N$ given by
\[\gamma(n)=|\{x\in X:x=s_1\dots s_nx_0\text{ for some }s_i\in S\}|.\]
$X$ has \emdef{polynomial growth} if $\gamma(n)\le n^D$ for some large
enough $D$, has \emdef{exponential growth} if $\gamma(n)\ge\lambda^n$
for some small enough $\lambda>1$, and has \emdef{intermediate growth}
in the remaining cases. The growth of $G$ is its growth under its
action on itself by left multiplication.

We describe now briefly five archetypical examples of groups.
%---------------------------------------------------------------
\subsection{The Groups $\Gg$ and $\tilde\Gg$}\label{subs:defG}
The first group was introduced by the second author in
1980~\cite{grigorchuk:burnside}; both groups act on $\tree_2$. Let $a$
be the automorphism permuting the top two branches of $\tree_2$, and
define recursively $b,c,d$ by $\phi(b)=(a,c)$, $\phi(c)=(a,d)$ and
$\phi(d)=(1,b)$.  Let $\Gg$ be the group generated by $\{a,b,c,d\}$.

Define also $\tilde b,\tilde c,\tilde d$ by $\phi(\tilde b)=(a,\tilde
c)$, $\phi(\tilde c)=(1,\tilde d)$ and $\phi(\tilde d)=(1,\tilde b)$,
and let $\tilde\Gg$ be the group generated by $\{a,\tilde b,\tilde
c,\tilde d\}$.  Clearly $\tilde\Gg$ contains $\Gg$ as the subgroup
$\langle a,\tilde b\tilde c,\tilde c\tilde d,\tilde d\tilde b\rangle$.

%----------------------------------------------------------------
\subsection{$GGS$ Groups}
Let $d$ be a prime number. Denote by $a$ the automorphism of
$\tree_d$ permuting cyclically the top $d$ branches. Fix a sequence
$\epsilon=(\epsilon_1,\dots,\epsilon_{d-1})\in(\Z/d)^{d-1}$. Define
recursively the automorphism $t_\epsilon$ of $\tree_d$, written
$t_\epsilon=(a^{\epsilon_1},\dots,a^{\epsilon_{d-1}},t_\epsilon)$, by
\[t_\epsilon(\sigma_1\sigma_2\dots\sigma_n)=\begin{cases}
  \sigma_1a^{\epsilon_{\sigma_1}}(\sigma_2\dots\sigma_n) &\text{ if }1\le\sigma_1\le d-1,\\
  \sigma_1t_\epsilon(\sigma_2\dots\sigma_n) &\text{ if }\sigma_1=d.
\end{cases}\]
Then $G_\epsilon$ is the subgroup of $\aut(\tree_d)$ generated by
$\{a,t_\epsilon\}$. It is called a \emdef{$GGS$
  group}~\cite{baumslag:cgt}.

The following results belong to folklore (see~\cite{grigorchuk:jibg}):
$G_\epsilon$ is an infinite group if an only if
$\epsilon\neq(0,\dots,0)$. It is a torsion group if and only if
$\sum\epsilon_i=0$. The only three infinite $GGS$ groups for $d=3$ are
as follows:

Let $r$ be the automorphism of $\tree_3$ defined recursively by
$\phi(r)=(a,1,r)$, and let $\Gamma$ be the subgroup of $\aut(\tree_3)$
generated by $\{a,r\}$. It was first considered by Narain Gupta and
Jacek Fabrykowski~\cite{fabrykowski-g:growth1}.

Define $s$ by $\phi(s)=(a,a,s)$, and let $\overline\Gamma$ be the
subgroup of $\aut(\tree_3)$ generated by $\{a,s\}$.

Define $t$ by $\phi(t)=(a,a^{-1},t)$, and let $\doverline\Gamma$ be
the subgroup of $\aut(\tree_3)$ generated by $\{a,t\}$; it was first
studied in the 80's by Narain Gupta and Said
Sidki~\cite{gupta-s:3group,gupta-s:infinitep}.

\begin{thm}
  The groups $\Gg,\tilde\Gg$ and all infinite $GGS$ groups with
  $d\ge3$ have intermediate growth. The Schreier graph $\sch(G,P,S)$
  corresponding to $G=\Gg,\tilde\Gg$ or any $GGS$ group has polynomial
  growth.
\end{thm}
In particular, $\sch(\Gg,P,S)$ and $\sch(\tilde\Gg,P,S)$ have linear
growth, while $\sch(\Gamma,P,S)$, $\sch(\overline\Gamma,P,S)$ and
$\sch(\doverline\Gamma,P,S)$ have polynomial growth of degree
$\log_2(3)$.

%----------------------------------------------------------------
\subsection{Substitutional Graphs}\label{subs:subs}
We give a self-contained description of the Schreier graphs
$\sch(G,P,S)$ for the examples above, in the form of
\emdef{substitutional rules}. In this subsection, all graphs shall
have a base point, and shall be edge-labelled; graph embeddings must
preserve the labelings.

\begin{defnn}
  A \emdef{substitutional rule} is a tuple $(U,R_1,\dots,R_n)$, where
  $U$ is a finite $d$-regular edge-labelled graph, called the
  \emdef{axiom}, and each $R_i$ is a rule of the form $X_i\to Y_i$,
  where $X_i$ and $Y_i$ are finite edge-labelled graphs. The graphs
  $X_i$ are required to have no common label. Furthermore, there is an
  inclusion, written $\iota_i$, of the vertices of $X_i$ in the
  vertices of $Y_i$; the degree of $\iota_i(x)$ is the same as the
  degree of $x$ for all $x\in V(X_i)$, and all vertices of $Y_i$ not
  in the image of $\iota_i$ have degree $d$.
\end{defnn}

Given a substitutional rule, one sets $\gf_0=U$ and constructs
iteratively $\gf_{n+1}$ from $\gf_n$ by listing all embeddings of all
$X_i$ in $\gf_n$ (noting that they are disjoint), and replacing them
by the corresponding $Y_i$. If the base point $*$ of $\gf_n$ is in a
graph $X_i$, the base point of $\gf_{n+1}$ will be $\iota_i(*)$.

Note that this expansion operation preserves the degree, so $\gf_n$ is
a $d$-regular finite graph for all $n$. We are interested in fixed
points of this iterative process, or equivalently in a converging
sequence of balls of increasing radius in the $\gf_n$, and call a
limit graph (which exists by~\cite{grigorchuk-z:infinite}) a
\emph{substitutional graph}.

\begin{thm}
  For the five examples $G$ described above, the Schreier graphs
  $\sch(G,P,S)$ are substitutional graphs.
\end{thm}

As an illustration, here is the substitutional rule for the group
$\Gamma$:
\begin{center}
  \begin{picture}(300,100)(0,-10)
    \put(0,0){\blue{\line(1,0){100}}}\put(50,-4){\blue{\msmash a}}
    \put(0,0){\blue{\line(15,26){50}}}\put(15,43){\blue{\msmash a}}
    \put(100,0){\blue{\line(-15,26){50}}}\put(85,43){\blue{\msmash a}}
    \put(0,-3){\footnotesize\msmash\rho}
    \put(100,-3){\footnotesize\msmash\sigma}
    \put(50,92){\footnotesize\msmash\tau}
    \put(100,40){\vector(1,0){40}}
    \put(130,0){\blue{\line(1,0){40}}}\put(150,7){\blue{\msmash a}}
    \put(130,0){\blue{\line(15,26){20}}}\put(143,17){\blue{\msmash a}}
    \put(170,0){\blue{\line(-15,26){20}}}\put(157,17){\blue{\msmash a}}
    \put(130,-3){\footnotesize\msmash{2\rho}}
    \put(144,37){\footnotesize\msmash{1\rho}}
    \put(168,-3){\footnotesize\msmash{0\rho}}
    \put(190,0){\blue{\line(1,0){40}}}\put(210,7){\blue{\msmash a}}
    \put(190,0){\blue{\line(15,26){20}}}\put(203,17){\blue{\msmash a}}
    \put(230,0){\blue{\line(-15,26){20}}}\put(217,17){\blue{\msmash a}}
    \put(195,-3){\footnotesize\msmash{1\sigma}}
    \put(217,37){\footnotesize\msmash{0\sigma}}
    \put(230,-3){\footnotesize\msmash{2\sigma}}
    \put(160,52){\blue{\line(1,0){40}}}\put(180,59){\blue{\msmash a}}
    \put(160,52){\blue{\line(15,26){20}}}\put(173,69){\blue{\msmash a}}
    \put(200,52){\blue{\line(-15,26){20}}}\put(187,69){\blue{\msmash a}}
    \put(158,57){\footnotesize\msmash{0\tau}}
    \put(180,92){\footnotesize\msmash{2\tau}}
    \put(205,57){\footnotesize\msmash{1\tau}}
    \put(170,0){\red{\line(-10,52){10}}}\put(157,37){\red{\msmash s}}
    \put(170,0){\red{\line(40,35){40}}}\put(185,7){\red{\msmash s}}
    \put(160,52){\red{\line(50,-17){50}}}\put(195,48){\red{\msmash s}}
    \put(150,35){\red{\spline(0,0)(10,17)(-10,17)(0,0)}}
    \put(190,0){\red{\spline(0,0)(-20,0)(-10,-17)(0,0)}}
    \put(200,52){\red{\spline(0,0)(20,0)(10,-17)(0,0)}}
    \put(300,70){\msmash{\text{axiom}}}
    \put(280,0){\blue{\line(1,0){40}}}\put(300,7){\blue{\msmash a}}
    \put(280,0){\blue{\line(15,26){20}}}\put(293,17){\blue{\msmash a}}
    \put(320,0){\blue{\line(-15,26){20}}}\put(307,17){\blue{\msmash a}}
    \put(283,-3){\footnotesize\msmash{2}}\put(272,7){\red{\msmash s}}
    \put(294,37){\footnotesize\msmash{1}}\put(328,7){\red{\msmash s}}
    \put(318,-3){\footnotesize\msmash{0}}\put(310,45){\red{\msmash s}}
    \put(300,35){\red{\spline(0,0)(10,17)(-10,17)(0,0)}}
    \put(280,0){\red{\spline(0,0)(-20,0)(-10,-17)(0,0)}}
    \put(320,0){\red{\spline(0,0)(20,0)(10,-17)(0,0)}}
  \end{picture}\\[2mm]
\end{center}

%%%%%%%%%%%%%%%%%%%%%%%%%%%%%%%%%%%%%%%%%%%%%%%%%%%%%%%%%%%%%%%%
\section{Amenability and Spectra}\label{sec:amen}
Let $G$ be a group acting on a set $X$. This action is
\emdef{amenable} in the sense of von Neumann~\cite{vneumann:masses} if
there exists a finitely additive measure $\mu$ on $X$, invariant under
the action of $G$, with $\mu(X)=1$.

A group $G$ is \emdef{amenable} if its action on itself by
left-multiplication is amenable.

We now state the main connection between the spectra of our
representations and dynamical systems. Recall $\pi$ is the
representation of $G$ on $L^2(\Sigma^\N)$ and $\pi_n$ is the
representation of $G$ on $L^2(\Sigma^n)$. Since $\pi_n$ contains
$\pi_{n-1}$ we write $\pi_n=\pi_{n-1}\oplus\pi_n^\perp$.
\begin{thm}
  Let $G$ be a group acting on a regular rooted tree, with $\pi$,
  $\pi_n$ and~$\pi^\perp_n$ be as above.
  \begin{enumerate}
  \item $\pi$ is a reducible representation of infinite dimension but
    splits as $\pi_0\oplus\bigoplus_{n\ge1}\pi_n^\perp$, so all of its
    irreducible components are finite-dimensional. Moreover
    $\spec(\pi)=\overline{\bigcup_{n\ge0}\spec(\pi_n)}$.
    
    However, if $G$ is weak branch
    (see~\cite{bartholdi-g:cras-parabolic}), then $\rho_{G/P}$ is
    irreducible.
  \item The representations $\pi_n$ and $\rho_{G/P_n}$ are equivalent,
    so their spectra coincide. The spectrum of $\rho_{G/P}$ is
    contained in $\overline{\cup_{n\ge0}\spec(\pi_n)}$, and therefore
    is contained in the spectrum of $\pi$.
    
    If moreover either $P$ or $G/P$ are amenable, these spectra
    coincide, and if $P$ is amenable, they are contained in the
    spectrum of $\rho_G$.
  \item $H_\pi$ has a pure-point spectrum, and its spectral radius
    $r(H_\pi)=s\in\R$ is an eigenvalue, while the spectral radius
    $r(H_{\rho_{G/P}})$ is not an eigenvalue of $H_{\rho_{G/P}}$.
    Therefore $H_{\rho_{G/P}}$ and $H_\pi$ are different operators
    having the same spectrum.
  \end{enumerate}
\end{thm}

%%%%%%%%%%%%%%%%%%%%%%%%%%%%%%%%%%%%%%%%%%%%%%%%%%%%%%%%%%%%%%%%
\section{Results on Spectra}
Since all the groups considered have intermediate growth, they are
amenable and $\spec\rho_{G/P}=\spec\pi$. We now describe explicitly
these spectra. For $\lambda\in\R$ let $J(\lambda)$ be the Julia set of
the quadratic map $z\mapsto z^2-\lambda$:
\[J(\lambda)=\overline{\left\{\sqrt{\lambda\pm\sqrt{\lambda\pm\sqrt{\lambda\pm\sqrt{\dots}}}}\right\}}.\]

\begin{center}
\begin{tabular}{r|c|p{6cm}}
Group & Spectrum of $M$ & Description\\[1pt] \hline
$\Gg$ & $[-\frac12,0]\cup[\frac12,1]$ &two intervals\rule{0pt}{12pt}\\
$\tilde\Gg$ & $[0,1]$ &nonnegative interval\\
$\Gamma$ & $\{1,\frac14\}\cup \frac14(1\pm J(6))\}$ &union of Cantor
set of null Lebesgue measure and set of isolated points\\
$\overline\Gamma$ &
$\{1,-\frac12,\frac14\}\cup\frac14\Big(1\pm\sqrt{\frac92\pm2J(\frac{45}{16})}\Big)$
&Cantor set of null Lebesgue measure\\
$\doverline\Gamma$ & & same as for $\overline\Gamma$\\
\end{tabular}
\end{center}

These computations imply the following results:
\begin{thm}
  \begin{enumerate}
  \item There are connected $4$-regular graphs of polynomial growth,
    which is are the Schreier graphs of groups of intermediate growth,
    and whose Markov operator's spectrum is any of the above sets.
  \item There are noncommutative dynamical systems generated by $3$
    (in the case of $\Gg$), $4$ (in the case of $\tilde\Gg$) or $2$
    transformations, whose spectrum is any of the above sets.
  \end{enumerate}
\end{thm}

These spectra are all computed using the same technique: consider the
representation $\pi_n$ of dimension $d^n$, given by $d^n\times
d^n$-permutation matrices $\pi_n(s)$, for all $s\in S$. These matrices
satisfy block identities, for instance for the group $\Gg$
\begin{xalignat*}{2}
  \pi_n(a)&=\MAT0{1_{d^{n-1}}}{1_{d^{n-1}}}0,&
  \pi_n(b)&=\MAT{\pi_{n-1}(a)}00{\pi_{n-1}(c)},\\
  \pi_n(c)&=\MAT{\pi_{n-1}(a)}00{\pi_{n-1}(d)},&
  \pi_n(d)&=\MAT{1_{d^{n-1}}}00{\pi_{n-1}(b)}.
\end{xalignat*}
Note that for our five examples $\pi_n(s)$ is expressed by blocks of
the form $0$, $1$ and $\pi_{n-1}(s')$.

Define now the polynomial $Q:\C^{S+1}\to\C$ by
\[Q_n(\{X_s\},\lambda)=\det\left(\sum_{s\in S}X_s\pi_n(s)-\lambda\right).\]
The spectrum of $\pi_n$ is
$\{\lambda\in\C|\,Q_n(\frac1{|S|},\dots,\frac1{|S|},\lambda)=0\}$.
Using the above block identities, it is possible to express $Q_n$
as rational expression over $Q_{n-1}$, and therefore to compute
inductively $Q_n$ and the spectrum of $\pi_n$. For the group $\Gg$,
for instance, we have for $n\ge2$, setting $\alpha=X_a$ and
$\beta=X_b=X_c=X_d$,
\[Q_n(\alpha,\beta,\lambda) = (3\beta^2+2\lambda\beta-\lambda^2)^{2^{n-2}}Q_{n-1}\left(\frac{2\alpha^2}{2\beta^2\lambda-\lambda^2},\lambda-\beta+\frac{(\lambda-\beta)\alpha^2}{3\beta^2+2\beta\lambda-\lambda^2}\right).\]

In all our cases the polynomials $Q_n$ can be explicitly computed;
again for the group $\Gg$, we have the
\begin{lem} Define
  \begin{xalignat*}{2}
    \Phi_0 &= \alpha+3\beta-\lambda,&
    \Phi_1 &= -\alpha+3\beta-\lambda,\\
    \Phi_2 &= -\alpha^2-3\beta^2-2\beta\lambda+\lambda^2,&
    \Phi_n &= \Phi_{n-1}^2 - 2(2\alpha)^{2^{n-2}}\text{ for }n\ge3.
  \end{xalignat*}
  Then for all $n\in\N$ we have $Q_n(\alpha,\beta,\lambda) =
  \Phi_0\Phi_1\cdots\Phi_n$.
\end{lem}

We show below the set of vanishing points of
$Q_6(X_a=X_{a^{-1}}=\alpha,X_r=X_{r^{-1}}=1,\lambda)$ for the group
$\Gamma$, and the Schreier graph
$\sch(\Gamma,P_6,\{a^{\pm1},r^{\pm1})$.
\begin{center}
\setlength\unitlength{1pt}
\begin{picture}(180,200)
\put(80,190){\makebox(0,0)[lb]{\smash{\normalsize $\lambda$}}}
\put(160,75){\makebox(0,0)[lb]{\smash{\normalsize $\alpha$}}}
\put(-60,-50){\epsfig{file=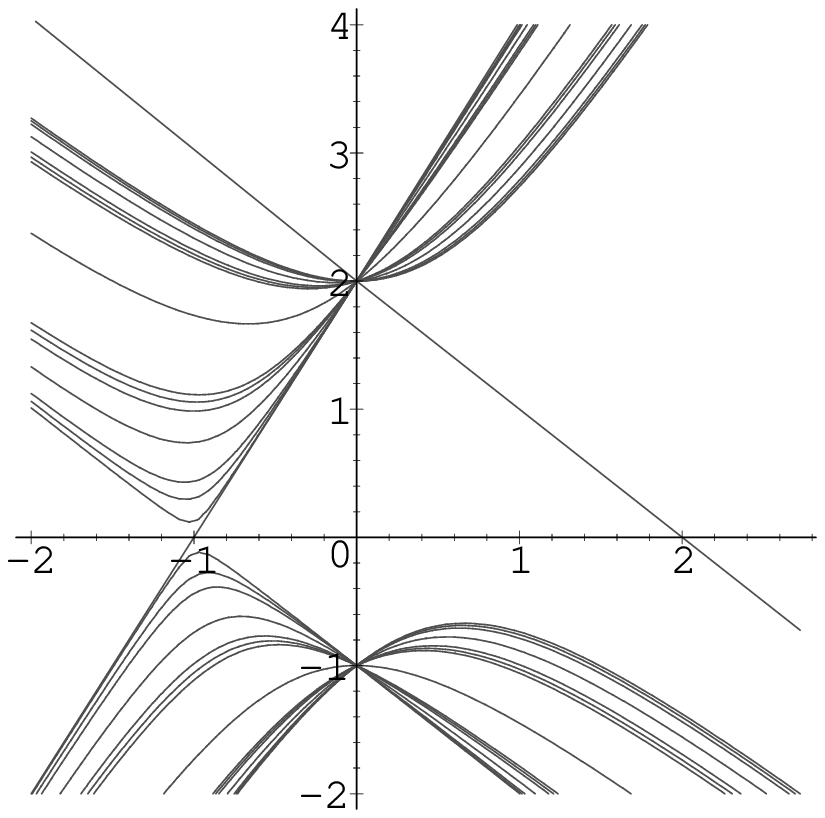,width=300pt}}
\end{picture}\qquad
\begin{picture}(180,200)
\put(-70,-65){\epsfig{file=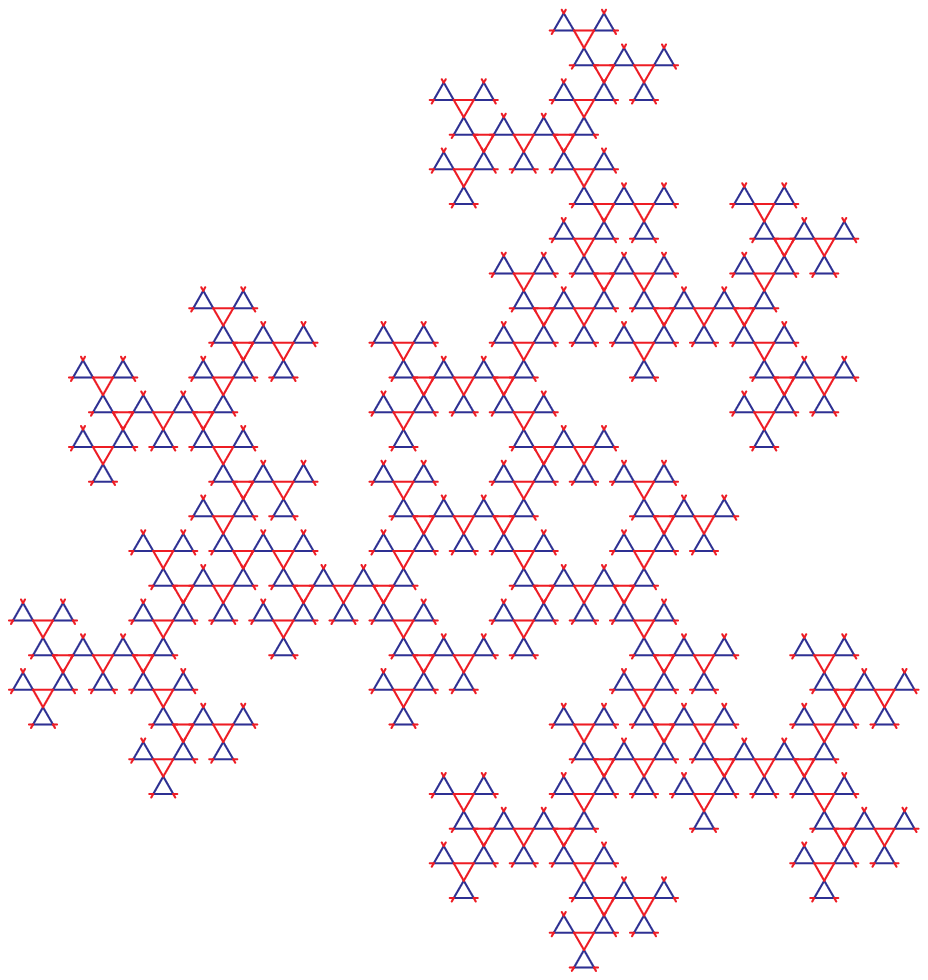,width=330pt}}
\end{picture}
\end{center}

\textit{\selectlanguage{francais}On montre ci-dessus l'ensemble des
  points d'annulation de
  $Q_6(X_a=X_{a^{-1}}=\alpha,X_r=X_{r^{-1}}=1,\lambda)$ pour le groupe
  $\Gamma$, et le graphe de Schreier
  $\sch(\Gamma,P_6,\{a^{\pm1},r^{\pm1})$.}

\bibliographystyle{amsplain}
\bibliography{mrabbrev,people,math,grigorchuk,bartholdi}
\end{document}